\documentclass[a4paper, 10pt]{article}
\usepackage[latin1]{inputenc}
\usepackage{amscd, amssymb}
\usepackage[mathscr]{eucal}
\usepackage{mathrsfs}
\usepackage{amsfonts}
\usepackage{amsmath}
\usepackage{amsthm}
\usepackage{latexsym}
\newcommand{\tr}{\mathrm{tr}|_{{\scriptscriptstyle M}}}
\newcommand{\C}{\mathbf{C}}
\newcommand{\N}{\mathbf{N}}
\newcommand{\CU}{\mathbf{C}[[U]]}
\newcommand{\Ct}{\mathbf{C}[[t_1,...,t_n]]}
\newcommand{\pf}{\textsl{Proof.}}
\newcommand{\A}{A\otimes A^o}
\newcommand{\epf}{\hspace{14,15 cm} $\Box$}

\newcommand{\Wr}{\C \Gamma \sharp \C\left\langle x,y\right\rangle/(xy-yx=\lambda)}

\newcommand{\Sm}{\C \Gamma \sharp \C\left\langle x,y\right\rangle}
\newcommand{\CSm}{\C \Gamma \sharp \C\left[ x,y \right]}
\newcommand{\G}{\mathbf{\Gamma}}
\newcommand{\g}{\gamma}
\newcommand{\h}{H_{t,k,c}(\G_N)}
\newcommand{\hh}{H_{1,k,c}(\G_N)}
\newcommand{\End}{\mathrm{End}}
\newcommand{\Ext}{\mathrm{Ext}}
\newcommand{\Hom}{\mathrm{Hom}}
\newcommand{\im}{\mathrm{Im}}
\newcommand{\Ker}{\mathrm{Ker}}
\newcommand{\rk}{\mathrm{rk}}
\newcommand{\ind}{\mathrm{Ind}}
\newcommand{\di}{\mathrm{dim}\,}
\newtheorem{defi}{Definition}[section]
\newtheorem{thm}[defi]{Theorem}
\newtheorem{lem}[defi]{Lemma}
\newtheorem{prop}[defi]{Proposition}
\newtheorem{coro}[defi]{Corollary}
\usepackage{latexsym}
\begin{document}

\title{Finite dimensional representations of symplectic
reflection algebras associated to wreath products}

\author{Pavel Etingof and Silvia Montarani}

\maketitle

\section{Introduction}
In this paper we study finite dimensional representations of
the wreath product symplectic reflection algebra
$H_{1,k,c}(\G_N)$ of rank $N$ attached to the group
$\G_N=S_N\ltimes \Gamma^N$ (\cite{EG}), where $\Gamma\subset
SL(2,\C)$ is a finite subgroup, and $(k,c)\in C(\mathcal{S})$, where
$C(\mathcal{S})$ is the space of (complex valued ) class
functions on the set $\mathcal{S}$ of symplectic reflections of
$\G_N$. 

In the rank $1$ case, there is no parameter $k$ and finite
dimensional representations of the wreath product algebra have
been classified in \cite{CBH} by Crawley-Boevey and Holland, by
establishing a Morita equivalence between the algebra
$H_{1,c}(\Gamma)$ and the deformed preprojective algebra
$\Pi_{\lambda(c)}(Q)$ attached to the (extended Dynkin) quiver $Q$
associated to $\Gamma$ via the McKay correspondence. 

We consider the higher rank case.  When $k=0$, we have 
$H_{1,k,c}(\G_N)=S_N\sharp H_{1,c}(\Gamma)^{\otimes N}$,
so the finite dimensional representations of $H_{1,k,c}(\G_N)$
are known. Using a cohomological approach, we investigate the
possibility of deforming some of these representations to values
of the parameters with $k\neq 0$. This allows us to produce the
first nontrivial examples of finite dimensional representations 
of $H_{1,k,c}(\G_N)$ for non-cyclic $\Gamma$ and $k\ne 0$. 

Specifically, we show
that if $W$ is an irreducible representation of $S_N$ 
whose Young diagram is a rectangle, and $Y$ an irreduible finite
dimensional representation of $H_{1,c}(\Gamma)$, then the
representation $M=W\otimes Y^{\otimes N}$ of 
$H_{1,0,c}(\G_N)$ can be deformed along a hyperplane in
$C({\mathcal S})$. On the other hand, if $\di Y=1$ and 
the Young diagram of $W$ is not a rectangle, such
a deformation does not exist.  
 
{\bf Acknowledgments.} The work of P.E. was partially supported
by the NSF grant DMS-9988796 and the CRDF grant RM1-2545-MO-03.

\section{Preliminaries}

\subsection{The wreath product construction}

In this subsection we will briefly recall the wreath product
construction. 
Let $L$
be a $2$-dimensional complex vector space with a symplectic form
$\omega_L$, and consider the space $V=L^{\oplus N}$,
endowed with the induced symplectic form
$\omega_V={\omega_L}^{\oplus N}$. Let $\Gamma$ be a finite
subgroup of $Sp(L)$, and let $S_N$ be the symmetric group acting
on $V$ by permuting the factors. The group
$\G_{N}:=S_N\ltimes\Gamma^N \subset Sp(V)$ acts naturally on
$V$. In the sequel we will write $\gamma_i\in \G_N$ for any
element $\gamma\in \Gamma$ seen as an element in the i-th factor
$\Gamma$ of $\G_N$. The symplectic reflections in $\G_N$ are the
elements $s$ such that $\rk(Id-s)|_V=2$. $\G_N$ acts by conjugation
on the set $\mathcal{S}$ of its symplectic reflections. It is
easy to see that there are symplectic reflections
of two types in $\G_N$:

\medskip

(S) the elements $s_{ij}\g_i {\g_j}^{-1}$ where $i,j\in [1,N]$, $s_{ij}$ is the transposition $(ij)\in S_N$, and $\g\in\Gamma$,

($\Gamma$) the elements $\g_i$, for $i\in[1,N]$ and $\g\in\Gamma\backslash\{1\}$.
\medskip

Elements of type (S) are all in the same conjugacy class, while
elements of type ($\Gamma$) form one conjugacy class for any
conjugacy class of $\g$ in $\Gamma$. 
Thus elements $f\in \C[\mathcal{S}]$ can be written as pairs
$(k,c)$, where $k$ is a number (the value of $f$ on elements of
type (S)), and $c$ is a conjugation invariant function
on $\Gamma\setminus \lbrace{1\rbrace}$ (encoding the values of
$f$ on elements of type ($\Gamma$)).
    
For any $s\in\mathcal{S}$ we write $\omega_s$ for
the bilinear form on $V$ that coincides with $\omega_V$ on
$\im(Id-s)$ and has $\Ker(Id-s)$ as radical. Denote by $TV$ the
tensor algebra of $V$.

\begin{defi} \label{dewre}
 For any $t\in \C$ and $f=(k,c)\in \C[\mathcal{S}]$, the symplectic reflection algebra $\h$ is  the quotient $$
(\G_N\sharp TV)/\left\langle [u,v]-\kappa(u,v)\right\rangle_{u,v\in V}
$$ 
where 
$$
\kappa:V\otimes V\longrightarrow\C[\G_N]: (u,v)\mapsto t\cdot \omega(u,v)\cdot 1 +\sum_{s\in \mathcal{S}}f_s\cdot \omega_s(u,v)\cdot s
$$
 with $f_s=f(s)$, and $\left\langle \dots\right\rangle$ is the two-sided ideal in 
the smash product $\G_N\sharp TV$ generated by the elements
$[u,v]-\kappa(u,v)$ for $u,v\in V$.

\end{defi} 

We will be interested in the case $t\neq 0$, and  it will be enough
to consider the case $t=1$ since $\h\cong H_{1,k/t,c/t}(\G_N)$ for any $t\neq 0$ (cf. \cite{EG}, page 14). We recall that the case $t=0$ is remarkably different and the corresponding representation theory has been studied in \cite{EG}, Section 3 and in \cite{GS}.

It is clear that choosing a symplectic basis $x$, $y$ for $L$ we
can consider $\Gamma$ as a subgroup of $SL(2,\C)$. We will denote by $x_i$, $y_i$ the corresponding vectors in the i-th $L$-factor of $V$.
Following \cite{GG} we will now  give a more explicit representation of the algebra $\hh$.
\begin{lem}([GG]) \label{pres}
The algebra $\hh$ is the quotient
of $\G_N\sharp TV$ by the following relations:
\begin{itemize}
\item[{\emph{(R1)}}] 
For any $i\in [1,N]$: $$[x_{i}, y_{i}]
=  1+ \frac{k}{2} \sum_{j\neq i}\sum_{\g\in\Gamma}
s_{ij}\g_{i}\g_{j}^{-1} + \sum_{\g\in\Gamma\backslash \{1\}}
c_{\g}\g_{i}\,.$$
\item[{\emph{(R2)}}] 
For any $u,v\in L$ and $i\neq j$:
$$[u_{i},v_{j}]= -\frac{k}{2} \sum_{\g\in\Gamma} \omega_{L}(\g u,v)
s_{ij}\g_{i}\g_{j}^{-1} \,.$$
\end{itemize}
\end{lem} 

\epf

In the case $N=1$, there is no parameter $k$ (there are no symplectic reflections of type (S)) and  

$$
H_{1,c}(\Gamma)=\Wr 
$$

where
$\lambda=\lambda(c)=1+\sum_{\g\in\Gamma\backslash\{1\}}c_{\g}\g\in
Z(\C[\Gamma])$ is the central element coresponding to $c$.

We end this subsection by recalling an important result 
that we will need in the sequel.
It is stated in \cite{EG} and is called the \textsl{Poincar\'e-Birkhoff-Witt (PBW-) property} for $\hh$.
Consider the increasing filtration on $TV\sharp \G_N$ obtained by
assigning degree zero to the elements of the group algebra
$\C[\G_N] $  and degree one to the vectors in $V$. 
This filtration induces a filtration on $\hh$. The following theorem holds:
\begin{thm}(PBW)\label{PBW}
The associated graded algebra to $\hh$ with respect to the above increasing filtration is $\G_N\sharp SV$, where $SV$ is the symmetric algebra of $V$.
\end{thm}

\epf

\subsection{Representations of $S_N$ with rectangular Young
diagram}
 
 We will use the following standard results from  representation
theory of the symmetric group. The proofs are well known, but we
recall them for reader's convenience. Denote by $\mathfrak{h}$ the
reflection representation of $S_N$. For a Young diagram $\mu$ we
denote by $\pi_{\mu}$ the corresponding irreducible
representation of $S_N$. 

\begin{lem}\label{sym}
\begin{center}
\begin{itemize}
\item[{(i)}] $\Hom_{S_N}(\mathfrak{h}\otimes \pi_{\mu},\pi_{\mu})=\C^{m-1}$, where $m$ is the number of corners of the Young diagram $\mu$. In particular  $\Hom_{S_N}(\mathfrak{h}\otimes \pi_{\mu},\pi_{\mu})=0$ if and only if $\mu$ is a rectangle.
\item[{(ii)}] The element $C=s_{12}+s_{13}+ \cdots +s_{1n}$ acts by a scalar in $\pi_{\mu}$ if and only if $\mu$ is a rectangle.
\end{itemize}
\end{center}
\end{lem}

\pf\ It is well known that $\pi_{\mu}|_{S_{N-1}}=
\sum \pi_{\mu -j}$, where the sum is taken over the corners of
$\mu$ and $\mu -j$ is the Young diagram obtained from $\mu$ by
cutting off the corner $j$. Since $\mathfrak{h}\oplus
\C=\ind_{S_{N-1}}^{S_N}\C$, the assertion {\emph(i)} follows from
the Frobenius reciprocity.
To prove {\emph(ii)}, observe that $C$ commutes with $S_{N-1}$,
so acts by a scalar on each $\pi_{\mu-j}$. Thus, if $\mu$ is a
rectangle, $C$ acts as a scalar (as we have only one summand), and
the ``if'' part of the statement is proved. To prove the ``only
if'' part, let  $Z_N$ be the sum of all transpositions in
$S_N$. $Z_N$ is a central element in the group algebra, and it is
known to act in $\pi_{\mu}$ by the scalar ${\bold c}(\mu)$, where
${\bold c}(\mu)$ is the content of $\mu$, i.e. the sum over all cells of
the signed distances from these cells to the diagonal. Now,
$C=Z_N-Z_{N-1}$, so it acts on $\pi_{\lambda-j}$ by the scalar
${\bold c}(j)$, the signed distance from the cell $j$ to the
diagonal. The numbers ${\bold c}(j)$ are clearly different for all
corners $j$, so if there are $2$ or more corners, then $C$ cannot
act by a scalar. This finishes the proof of (ii).

\epf

\section{The main theorem}

Let $Y$ be an irreducible representation 
of the algebra $H_{1,c}(\Gamma)$ for some $c$
\footnote{Such representations exist only for special $c$, as for
generic $c$ the algebra $H_{1,c}(\Gamma)$ is simple; see \cite{CBH}.}. 
Let $W$ be an irreducible representation of $S_N$.
Since the algebra $H_{1,0,c}(\bold \Gamma_N)$
is naturally isomorphic to $S_N\sharp H_{1,c}(\Gamma)^{\otimes
N}$, there is a natural action of $H_{1,0,c}(\bold\Gamma_N)$ 
on the vector space $M:=W\otimes Y^{\otimes N}$. 
Namely, each copy of $H_{1,c}(\Gamma)$ acts in the corresponding
copy of $Y$, while $S_N$ acts in $W$ and simultaneously permutes
the factors in the product $Y^{\otimes N}$. We will denote this
representation by $M_c$. 
The main theorem tells us 
when such a representation can deformed to nonzero values
of $k$. 

Assume that the Young diagram of $W$ is a rectangle 
of height $l$ and width $m=N/l$ (the trivial representation
corresponds to the horizontal strip of height $1$). 

Let ${\mathcal H}_{Y,m,l}$ be the hyperplane in $C(\mathcal S)$ 
consisting of all pairs $(k,c)$ satisfying the equation 
\begin{equation}\label{hyp}
\di Y+\frac{k}{2}|\Gamma|(m-l)+\sum_{\gamma\in
\Gamma\setminus \lbrace{1\rbrace}}c_\gamma\chi_Y(\gamma),
\end{equation}
where $\chi_Y$ is the character of $Y$. 

Let $X=X(Y,m,l)$ be the moduli space of irreducible representations
of $H_{1,k,c}(\bold\Gamma_N)$  isomorphic to $M$ 
as $\G_N$-modules (where $(k,c)$ are allowed to vary). 
This is a quasi-affine algebraic variety:
it is the quotient of the quasi-affine variety
$\widetilde{X}(Y,m,l)$ of extensions of 
the $\bold\Gamma_N$-module $M$ 
to an irreducible $H_{1,k,c}(\bold\Gamma_N)$-module 
by a free action of the reductive group $G$ of basis changes in $M$
compatible with $\G_N$ modulo scalars. Let $f: X\to
C({\mathcal S})$ 
be the morphism which sends a representation to the
corresponding values of $(k,c)$. 

The main result of this paper is the following theorem. 

\begin{thm}\label{mainth}
(i) For any $c_0$ the representation $M_{c_0}$ of $H_{1,0,c_0}(\bold \Gamma_N)$ 
can be formally deformed to a representation of $H_{1,k,c}(\bold
\Gamma_N)$ along the hyperplane ${\mathcal H}_{Y,m,l}$, but not
in other directions. This deformation is unique. 

(ii) The morphism $f$ maps $X$ to ${\mathcal H}_{Y,m,l}$ 
and is etale at $M_{c_0}$ for all $c_0$. Its restriction to the formal
neighborhood of $M_{c_0}$ is the deformation from (i).  

(iii) There exists a nonempty Zariski open subset 
$\mathcal U$ of ${\mathcal H}_{Y,m,l}$ such that for $(k,c)\in
\mathcal U$, the algebra $H_{1,k,c}(\bold\Gamma_N)$ admits a 
finite dimensional irreducible representation isomorphic to $M$ as a 
$\bold \Gamma_N$-module.  
\end{thm}

The proof of this theorem occupies the remaining sections of the paper.

{\bf Remark.} In the case of cyclic $\Gamma$
and trivial $W$ Theorem \ref{mainth} was proved in \cite{CE}. 
In this case, the deformation of the representation $M$ 
can be constructed explicitly.

 We expect that the condition that the Young diagram 
of $W$ is a rectangle is essential to obtain the deformation 
of Theorem \ref{mainth} (i). For example, this is the case if $Y$
is 1-dimensional. This follows from the following more general
statement. 

\begin{prop} Let $W$ be an irreducible 
$S_N$-module. If $W$ extends to a representation of
$H_{1,k,c}(\bold\Gamma_N)$ for some $(k,c)$ with $k\ne 0$, 
then the Young diagram of $W$ is a rectangle.
\end{prop}

\pf\  Suppose that $W$ extends to a representation of 
$H_{1,k,c}(\bold\Gamma_N)$. Such an extension is, first of all, 
an extension of $W$ to a 
representation of the wreath product group $S_N\ltimes \Gamma^N$. 
This can only be done by making $\Gamma^N$ act by an
$S_N$-invariant character $\xi$, i.e.,
$\xi(\gamma_1,...,\gamma_N)=\chi(\gamma_1)...\chi(\gamma_N)$,
where $\chi:\Gamma\to \C^*$ is a character. 
But in this case $\Gamma^N$ acts trivially on ${\rm End}_\C(W)$,
and hence $x_i, y_i$ must act by  $0$ on $W$ for each
$i=1,...,N$. 
So, denoting by $\rho$ the possible extended representation, we obtain from relation {\emph{(R1)}} for $i=1$:
$$ 
\rho\left(s_{12}+...+s_{1N}\right)=-2\,\frac{1+\sum_{\g\in\Gamma\backslash\{1\}}c_{\g}\chi(\gamma)}{k\,|\Gamma|}
$$
i.e. $C=s_{12}+...+s_{1N}$ acts by a constant on $W$. Now applying
Lemma \ref{sym}, part \textsl{(ii)}, 
we get that $W$ must correspond to a rectangular Young diagram.

\epf
\bigskip 

\section{Proof of Theorem \ref{mainth}}
\subsection{Deformation theory.}

In this section we recall deformation theory of representations
of algebras. This theory is well known, but 
we give the details for reader's convenience. 

Let $A$ be an associative algebra over $\C$. In what follows, for each $A$-bimodule $E$, we write $H^n(A,E)$ for  the n-th Hochschild cohomology group of $A$ with coefficients in $E$. We recall that $H^n(A,E)$ is defined to be the i-th cohomology group of the Hochschild complex :
$$
0\longrightarrow C^0(A,E)\stackrel{d}{\longrightarrow} \cdots \stackrel{d}{\longrightarrow} C^n(A,E)\stackrel{d}{\longrightarrow} C^{n+1}(A,E)\stackrel{d}{\longrightarrow} \cdots
$$
where $C^n(A,E)=\Hom_{\C}(A^{\otimes n},E)$ is the space of n-linear maps from $A^n$ to $E$, and the differential $d$ is defined as follows:
\begin{eqnarray*}
(d\varphi)(a_1,\cdots,a_{n+1}):&=& a_1\varphi(a_2,\cdots,a_{n+1})\\
                               &+& \sum_{i=1}^n (-1)^i \varphi(a_1,\cdots,a_{i-1},a_i\,a_{i+1},a_{i+2},\cdots,a_{n+1})\\
                               &-& (-1)^n \varphi(a_1,\cdots,a_n)\,a_{n+1}.
\end{eqnarray*}
We remark that    $H^i(A,E)$  coincides with the vector space $Ext^i_{\A}(A, E)$, where $A^o$ is the opposite algebra of $A$.

Let $A_U$ be a flat formal deformation of $A$ over the formal
neighborhood of zero in a
finite dimensional vector space $U$ with coordinates
$t_1,...,t_n$. This means that $A_U$ is an algebra over $\CU =
\Ct$ which is \textsl{topologically free} as a $\CU$-module (i.e., $A_U$  is isomorphic as a $\CU$-module to $A[[U]]$), together with a fixed isomorphism of algebras $A_U/JA_U\cong A$, where $J$ is the maximal ideal in $\CU$.  
Given such a deformation, we have a natural linear map\  $\phi:U \longrightarrow H^2(A, A)$.

Explicitly, we can think of $A_U$ as $A[[t_1,...,t_n]]$ equipped
with a new $\Ct$-linear (and continuous) associative product defined by: 
$$
a \ast b =\sum_{p_1,...,p_n}c_{p_1,...,p_n}(a,b)\prod_j t_j^{p_j} \qquad a,b \in A  
$$
where $c_{p_1,...,p_n}:A\times A\longrightarrow A$ are
$\C$-bilinear functions and $c_{0,...,0}(a,b)=ab$, for any
$a,b\in A$. 

Imposing the associativity condition on $\ast$, one obtains that $c_{0,...,1_j,...0}$ must be Hochschild $2$-cocycles for each $j$.  
The map $\phi$ is given by the assignment
$(t_1,\cdots,t_N)\rightarrow \sum_jt_j\,[c_{0,...,1_j,...0}]$ for
any $(t_1,\cdots,t_n)\in U$, where $[C]$ stands for the
cohomology class of a cocycle $C$. 

Now let $M$ be a representation of $A$. In general it does not
deform to a representation of $A_U$. However we have the
following standard proposition. Let $\eta: U\to H^2(A, \End M)$ be the composition of $\phi$ with the natural map $\psi:H^2(A,A)\longrightarrow H^2(A, \End\,M)$. 

\begin{prop} \label{defo}
Assume that $\eta$ is surjective with kernel $K$, and $H^1(A,\End M)=0$. Then: 

(i) There
exists a unique smooth formal subscheme $S$ of the formal neighborhood
of the origin in 
$U$, with tangent space $K$ at the origin, 
such that  $M$ deforms to a representation of the algebra 
$A_S:=A_U \hat\otimes_{\CU}\C[S]$ (where $\hat\otimes$ is the
completed tensor product).

(ii) The deformation of $M$ over $S$ is unique.
\end{prop}
 
 \pf \ Let us realize $A_U$ explicitly  as $A[[t_1,... ,t_n]]$
equipped with a product $\ast$ as above.
We may assume that $K$ is the space of all vectors $(t_1,...,t_n)$
such that $t_{m+1}=...=t_n=0$. 

Let $D$ be the formal neighborhood 
of the origin in $K$, with coordinates $h_1=t_1,...,h_m=t_m$. Let  
$\theta: D\to U$ be a map given 
by the formula $\theta(h_1,...,h_m)=(t_1,...,t_n)$, 
where $t_i=h_i$ for $i\le m$, and 
$$
t_k=\sum_{p_1,...,p_m}t_{k,p_1,...,p_m}h_1^{p_1}...h_m^{p_m}, k>m,
$$
where $t_{k,p_1,...,p_m}\in \C$. 
More briefly, we can write $t_k=\sum_P t_{kP}h^P$, where 
$P$ is a multi-index. 
We will use the notation $|P|$ 
for the sum of indices in a multi-index $P$.
For brevity we also let $e_j$ to be the multi-index
$(0,...,1_j,...,0)$. 

We claim that there exist unique 
formal functions $t_k=t_k(h)$, $k>m$, for which we can deform $M$ over $D$. 
Indeed, such a deformation would be defined by a series 
$$
\widetilde{r}(a)=\sum_{P}r_{P}(a)h^P,
$$
where $r_0(a)=r(a)$, and $r$ is the homomorphism giving the
representation $M$. The condition that $\widetilde{r}$ is a representation
gives, for each $P$, 
\begin{equation}\label{drp}
d\,r_P=\sum_j t_{jP}r(c_{e_j})+C_P,
\end{equation}
where for $j\le m$, $t_{jP}=1$ if $P=e_j$ and zero otherwise,  
and $C_P$ is a $2$-cocycle
whose expression may involve $r_Q$ and $t_{kQ}$ 
{\bf only} with $|Q|<|P|$. Since the map $\eta$ is surjective,
there are (unique) $t_{d+1,P},...,t_{nP}$ for which the right hand side 
is a coboundary. For such $t_{d+1,P},...,t_{nP}$
(and only for them), we can solve (\ref{drp}) for $r_P$.

This shows the existence of the functions $t_j(h)$, $j>m$, such
that the deformation of $M$ over $D$ is possible. To show the
uniqueness of these functions, let $t_j$ and $t_j'$ be two sets
of functions for which the deformation exists. 
Let $r_P,r'_P$ be the coeffients of the corresponding
representations $\widetilde{r},\widetilde{r}'$. Let $N$ be the
maximal number such that $t_{jP}=t'_{jP}$ for $|P|<N$. 
Since $H^1(A,\End M)=0$, the solution $r_P$ of (\ref{drp})
is unique up to adding a coboundary. Thus
we can use changes of basis in $M$ 
to modify $\widetilde{r}$ so that $r_P=r'_P$ for $|P|<N$ (note that this
does not affect $t_j$). 
Then for any $Q$ with $|Q|=N$, $C_Q(\widetilde{r})=C_Q(\widetilde{r'})$, and hence 
$t_{jQ}=t_{jQ}'$. This contradicts the maximality of $N$. 

Thus, we have shown that the functions $t_j$ exist and are
unique; they define a parametrization of the desired subscheme
$S$ by $D$. 
Our proof also implies that the deformation of $M$ over $S$ is
unique, so we are done. 

\epf

\subsection{Homological properties of the algebra $H_{1,c}(\Gamma)$.}
We recall the following definition (see \cite{VB1,VB2,EO}):

\begin{defi}\label{dclass}
An algebra $A$ is defined to be in the class $VB(d)$ if it is of finite Hochschild dimension (i.e. there exists $n\in\N$ s.t. $H^i(A,E)=0$ for any $i>n$ and any $A$-bimodule $E$) and $H^{\ast}(A,A\otimes A^o)$ is concentrated in degree $d$, where it equals $A$ as an $A$-bimodule.
\end{defi}

The meaning of this definition is clarified by the following
result by Van den Bergh (\cite{VB1,VB2}).

\begin{thm}\label{vdb}
If $A\in VB(d)$ then for any $A$-bimodule $E$, 
the Hochschild homology $H_i(A,E)$ is naturally isomorphic to the
Hochschild cohomology $H^{d-i}(A,E)$. 
\end{thm}

\epf

Now let $B=H_{1,c}(\Gamma)$. 

\begin{prop}\label{vander} 
The algebra $B$ belongs to the class $VB(2)$.
\end{prop}
 
\pf\ 
If $\Gamma=\lbrace{1\rbrace}$, the statement is well known (\cite{VB1,VB2}; see
also \cite{EO}). Let us consider the case $\Gamma\ne \lbrace{1\rbrace}$. 
 We have to show that $B$ has finite Hochschild dimension and that:
 
 $$
 H^i(B,B\otimes B^o)=0 \qquad \mbox{for\ } i \neq 2
 $$ 
 
 $$
 H^2(B,B\otimes B^o)\cong B \qquad \mbox{as\ } B-\mbox{bimodules}. 
 $$ 
The algebra $\Sm$ has a natural increasing filtration obtained
by putting $x$, $y$ in degree $1$ and the elements of $ \Gamma$
in degree $0$. This filtration clearly induces a filtration on
$B$: $B=\cup_{n\ge 0}F^nB$, 
and the associated graded algebra is $B_0=grB=\CSm$ (by the PBW
theorem), 
which has Hochschild dimension 2.
So by a deformation argument we have that $B$ has finite
Hochschild dimension (equal to $2$) and $H^i(B,B\otimes B^o)=0$
for $i\neq 2$, as this is true for $B_0$ (which is easily checked
since $B_0$ is a semidirect product 
of a finite group with a polynomial algebra). 

It remains to show the
$B$-bimodule $E:=H^2(B,B\otimes B^o)$ is isomorphic to $B$. 
Using again a deformation argument (cf. \cite{VB1}), 
we can see that $E$ is
invertible and free as a right and left $B$-module, because this is
true for $B_0$. So $E=B\phi$ 
where $\phi$ is an automorphism of $B$ such that $gr\phi =1$.
We will now show that $\phi=1$, which will conclude the proof. 

Define a linear map $\xi: B_0\to B_0$ as follows: 
if $z\in B_0$ is a homogeneous element of degree $n$, 
and $\widetilde z$ is its lifting to $B$, then 
$\xi(z)$ is defined to be the 
projection of the element $\phi(\widetilde z)-\widetilde z$
(which has filtration degree $n-1$) to $B_0[n-1]$.  
It is easy to check that $\xi$ is well defined (i.e., independent
on the choice of the lifting), and is a derivation of $B_0$ of
degree $-1$. 

Our job is to show that $\xi=0$. This would imply that $\phi=1$, 
since $B$ is generated by $F^1B$. 

It is clear that any homogeneous inner derivation of $B_0$
has nonnegative degree. Hence, it suffices to show that 
the degree $-1$ part of $H^1(B_0,B_0)$ is zero. 
But it is easy to compute using Koszul complexes that 
$H^1(B_0,B_0)={\rm Vect}(L)^\Gamma$, the space 
of $\Gamma$-invariant vector fields on $L$. 
In particular, vector fields of degree $-1$ are those with constant
coefficients. But such a vector field cannot be
$\Gamma$-invariant unless it is zero, since the space 
$L$ has no nonzero vectors fixed by $\Gamma$. 
Thus, $\xi=0$ and we are done. 

\epf

\begin{coro}\label{endoo}
$H^2(B,\End\,Y)=H_0(B, \End\,Y)=\C$.
\end{coro}
\pf \ We apply  Theorem \ref{vdb}, to obtain the first
identity. Furthermore, $H_0(B, \End\,Y)=\End\,Y/[B, \End\,Y]=\C$ as $Y$
is irreducible, 
so the second identity follows.

\epf 

\begin{prop}\label{ext1}
$H^1(B,\End Y)=0$.
\end{prop}

\pf\ 
We have $H^1(B, \End Y)=
\Ext^1_{B\otimes B^o}(B, \End\,Y)=\Ext^1_B(Y,Y)$. But it is known 
(\cite{CBH}, Corollary $7.6$) 
that $B$ contains only one minimal ideal $J$ among all the
nonzero ideals, and $\Ext^1_{B/J}(Y',Y')=0$ 
for any irreducible module $Y'$ over 
the (finite dimensional) quotient algebra $B/J$.
Since any finite dimensional $B$-module must factor through
$B/J$, we get $\Ext^1_B(Y,Y)=0$, as desired. 

\epf

\subsection{Homological properties of
$A=H_{1,0,c}(\bold\Gamma_N)$.} 

We now let $A$ denote the algebra $H_{1,0,c_0}(\bold\Gamma_N)$. 
The algebra $A$ has a flat deformation over 
$U=C({\mathcal S})$, which is given by the algebra 
$H_{1,k,c_0+c'}(\bold\Gamma_N)$. The fact that this deformation is
flat follows from Theorem \ref{PBW}. 

\begin{prop}\label{endot}
If the Young diagram of $W$ is a rectangle then
$$
H^2(A,\End\,M)=H^2(B, \End\,Y)=\C.
$$
\end{prop}

\pf\  The second equality follows from Corollory \ref{endoo}.
Let us prove the first equality. 
We have: 
$$
H^{\ast}(A,\End\,M)=\Ext^{\ast}_{A\otimes A^o}(A,\End\,M)=
$$
$$
=\Ext^{\ast}_{S_N\sharp B^{\otimes N}\otimes S_N\sharp {B^o}^{\otimes N}}(S_N\sharp B^{\otimes N},\End\,W\otimes{\End\,Y}^{\otimes N})=
$$
$$
=\Ext^{\ast}_{S_N\times S_N\sharp(B^{\otimes N}\otimes {B^o}^{\otimes N})}(S_N\sharp B^{\otimes N},\End\,W\otimes {\End\,Y}^{\otimes N}).
$$
Now, the  $S_N\times S_N\sharp(B^{\otimes N}\otimes {B^o}^{\otimes
N})$-module $S_N\sharp B^{\otimes N}$ is induced from the module $B^{\otimes N}$ over the subalgebra $S_N\sharp B^{\otimes N}\otimes {B^o}^{\otimes N}$, in which $S_N$ acts simultaneously  permuting the factors of $B^{\otimes N}$ and ${B^o}^{\otimes N}$ (note that $S_N\sharp (B^{\otimes N}\otimes {B^o}^{\otimes N})$ is indeed a subalgebra of $S_N\times S_N\sharp(B^{\otimes N}\otimes {B^o}^{\otimes N})$ as it can be identified with the subalgebra $D\sharp(B^{\otimes N}\otimes {B^o}^{\otimes N})$ where $D=\left\{(\sigma,\sigma),\; \sigma\in S_N\right\}\subset S_N\times S_N$). Applying the Shapiro Lemma, we get:
$$
 \Ext^{\ast}_{S_N\times S_N\sharp(B^{\otimes N}\otimes {B^o}^{\otimes N})}(S_N\sharp B^{\otimes N},\End\,W\otimes\End\,Y^{\otimes N})=
$$
$$
=\Ext^{\ast}_{S_N\sharp(B^{\otimes N}\otimes {B^o}^{\otimes N})}(B^{\otimes N},\End\, W\otimes \End\,Y^{\otimes N})=
$$
$$
={\left(\Ext^{\ast}_{B^{\otimes N}\otimes {B^o}^{\otimes N}}(B^{\otimes N},\End\, W\otimes \End\,Y^{\otimes N})\right)}^{S_N}.
$$
But since $B^{\otimes N}\otimes {B^o}^{\otimes N}$ does not act on $\End\,W$, the latter module equals:
$$
{\left(\Ext^{\ast}_{B^{\otimes N}\otimes {B^o}^{\otimes N}}( B^{\otimes N},{\End\,Y}^{\otimes N})\otimes \End\,W\right)}^{S_N}.
$$

  Using Proposition \ref{ext1} and the K\"unneth formula in degree $2$,  we get that as an $S_N$-module, $\Ext^2 _{B^{\otimes N}\otimes {B^o}^{\otimes N}}(B^{\otimes N},{\End\,Y}^{\otimes N})=\Ext^2_{B\otimes B^o}(B, \End\,Y)\otimes \C^N$
where $S_N$ acts only on $\C^N$ permuting the factors.
But as an $S_N$-module, $\C^N=\C\oplus\mathfrak{h}$, where $\C$ is the trivial representation.
As a result  we get:
$$
\Ext^2_{\A}(A,\End\,M)= \Ext^2_{B\otimes B^o}(B, \End\,Y)\otimes{\left(\C^N\otimes \End(W)\right)}^{S_N}=
$$
$$
=\Ext^2_{B\otimes B^o}(B, \End\,Y)\otimes{\left(\C\otimes \End(W)\oplus \mathfrak{h}\otimes \End\,W\right)}^{S_N}=
$$
$$
=\Ext^2_{B\otimes B^o}(B, \End\,Y)\otimes \left(\Hom_{S_N}(W,W)\oplus \Hom_{S_N}(\mathfrak{h}\otimes W, W)\right)=
$$
$$
=\Ext^2_{B\otimes B^o}(B, \End\,Y)
$$
as $\Hom_{S_N}(\mathfrak{h}\otimes W,W)=0$ by Lemma \ref{sym} part \emph{(i)}.

\epf

\begin{coro}\label{surj}
The map $\eta:U\longrightarrow H^2(A,\End\,M)$ is surjective.
\end{coro}
\pf\  
Let $U_0\subset U$ be the subspace of vectors $(0,c')$. 
It is sufficient to show that the restriction
of 
$\eta$ to $U_0$ is surjective. 
But this restriction is a composition of three natural maps:
$$
U_0\to H^2(B,B)\to H^2(A,A)\to H^2(A,\End M).
$$
Here the first map $\eta_0: U_0\to H^2(B,B)$
is induced by the deformation of $B$ along $U_0$, 
the second map $\xi: H^2(B,B)\to H^2(A,A)$ 
comes from the K\"unneth formula, 
and the third map $\psi: H^2(A,A)\to H^2(A,\End M)$ 
is induced by the homomorphism $A\to \End M$. 

Now, by Proposition \ref{endot}, 
the map $\psi\circ \xi$ coincides with the 
map $\psi_0: H^2(B,B)\to H^2(B,\End Y)$ induced 
by the homomorphism $B\to \End Y$. 
We claim that this map is surjective. 
Indeed, since by Proposition \ref{vander}, $B$ is in VB(2), by Theorem \ref{vdb} 
there is a natural identification 
of $H^2(B,E)$ with $H_0(B,E)$ for any $B$-bimodule $E$; hence 
$\psi_0$ can be viewed as the natural map 
$\psi_0: H_0(B,B)\to H_0(B,\End Y)$.  
But $H_0(B,E)=E/[B,E]$ for any $B$-bimodule $E$. 
Hence, $\psi_0$ can be viewed as the natural map 
$$
\psi_0: B/[B,B]\longrightarrow \End\,Y/[B,\End\,Y].
$$
This map is clearly nonzero: the representation
$Y$ is irreducible, and hence the map $B\to \End Y$ is
surjective. Thus $\psi_0$ is surjective, as claimed
(as the space $\End\,Y/[B,\End\,Y]$ is 1-dimensional).

Let $K$ be the kernel of $\psi_0$. 
It remains to show that the map $\eta_0$ does not land in $K$. 
To show this, recall that by Proposition \ref{defo},
the representation $Y$ of $B$ can be
deformed along $K$. Thus it remains to show that 
$Y$ does not admit a first order deformation along the entire $U_0$. 
But this follows easily by computing the trace of both sides 
of the commutation relation $xy-yx=\lambda$ in a deformation of
$Y$. We are done. 

\epf

\begin{prop}\label{h1}
$H^1(A,\End M)=0$. 
\end{prop}

\pf\ 
Arguing as in the proof of Proposition 
\ref{endot}, we get that $H^1(A, \End M)=H^1(B,\End Y)$, which is
zero. This proves the proposition. 

\epf

We have thus proved the following result. 

\begin{prop} \label{main}
If the Young diagram corresponding 
to $W$ is a rectangle, then there exists 
a unique smooth codimension one formal subscheme $S$
of the formal neighborhood of the origin in $U$
such that the representation $M=W\otimes Y^{\otimes N}$ of 
$H_{1,0,c_0}(\bold\Gamma_N)$ formally deforms 
to a representation of $H_{1,k,c_0+c'}(\bold\Gamma_N)$ 
along $S$ (i.e., abusing the language, for $(k,c')\in S$). 
Furthermore, the deformation of $M$ over $S$ is unique.
\end{prop}

\pf\ Corollary \ref{surj} and Proposition \ref{h1} show 
that our case satisfies 
all the hypothesis of Proposition \ref{defo}. Moreover, from  $H^2(A, \End\,M)=\C$ we deduce $\di\Ker\,\eta=\di U-1$, and the Proposition follows.

\epf

\subsection{The trace condition and the proof of Theorem \ref{mainth}}
 
Now we would like to find the subscheme $S$ of Proposition \ref{main}. 
For this we take the trace in $M$ of the 
commutation relation (\emph{R1}), and obtain a necessary condition on
the parameters $(k, c)$ for the algebra
$H_{1,k,c}(\bold\Gamma_N)$
to admit a representation 
isomorphic to $M$ as a $\Gamma_N$-module:
\begin{itemize}
\item[{\emph{(TR)}}] 
For any $i\in [1,n]$: $$0
=   \di M+ \frac{k}{2} \sum_{j\neq i}\sum_{\g\in\Gamma}
\tr(s_{ij}\g_{i}\g_{j}^{-1}) + \sum_{\g\in\Gamma\backslash\{1\}}
c_{\g}\tr(\g_{i})\,.$$
\end{itemize}

This relation can be easily rewritten in terms of the the characters $\,\chi_{{\scriptscriptstyle Y}}\,$ of $Y$ as a representation of  $\Gamma$ and $\,\psi_{{\scriptscriptstyle W}}$ of $W$ as a representation of $S_N$.
Indeed, one can check:
\begin{equation}\label{trf}
\tr(\g_i)={\di W}\,{\di Y}^{N-1}\,\chi_{\scriptscriptstyle Y}(\g)
\end{equation}
\begin{equation}\label{trfi}
\tr(s_{ij}\,\g_i\,{\g_j}^{-1})=\psi_{{\scriptscriptstyle W}}(s_{ij})\,{\di Y}^{N-1}
\end{equation}

Namely, (\ref{trf}) is an easy consequence of the fact that the
group $\Gamma^{\times N}\subset\G_N$ acts only on $Y^{\otimes N}$
with character $\chi_{\scriptscriptstyle Y}^{\otimes N}$, and $\g_i$ is by definition the element $(1,\dots ,\stackrel{i}{\g}, \dots, 1)\in \Gamma^{\times N}$.
To obtain (\ref{trfi}), we observe that $s_{ij}\,\g_i\,{\g_j}^{-1}$ is conjugate in $\G_N$ to $s_{ij}$ and that the character of $S_N$ on $M$ is simply the product of the characters on $W$ and $Y^{\otimes N}$. An easy computation gives $\mathrm{tr}|_{\scriptscriptstyle{Y^{\otimes N}}}s_{ij}={\di Y}^{N-1}$, hence the formula.

We now recall that, for any transposition $\sigma\in S_N$, $\psi|_{\scriptscriptstyle{W}}(\sigma)=\frac{\di W}{N\,(N-1)/2} {\bold c}(\mu)$, where ${\bold c}(\mu)$ is the content of the Young diagram $\mu$ attached to $W$. In particular, if $\mu$ is a rectangular diagram of size $l\times m$ with $lm=N$, it can be easily computed that: 
$$
{\bold c}(\mu)=\frac{N\,(m-l)}{2},
$$
so we have:
\begin{equation}\label{trsi}
\tr (s_{ij}\g_i\,{\g_j}^{-1})=\frac{(m-l)\,\di W}{N-1} {\di Y}^{N-1}
\end{equation}
Finally, substituting (\ref{trf}), (\ref{trsi}) in {\emph{(TR)}} and dividing the relation by ${\di Y}^{N-1}\,\di W$, we obtain:
\begin{itemize}
\item[{\emph{(TR')}}] 
If the Young diagram of $W$ is of size $l\times m$:
$$0
= \di Y+ \frac{k}{2}\;|\Gamma|(m-l)+ \sum_{\g\in\Gamma\backslash \{1\}}
c_{\g}\chi_{{\scriptscriptstyle Y}}(\g).$$ 
\end{itemize}
The condition {\emph{(TR')}} defines exactly
the hyperplane ${\mathcal H}_{Y,m,l}$. 

Thus we have shown that $(0,c_0)+S\subset {\mathcal H}_{Y,m,l}$. 
But $S$ and ${\mathcal H}_{Y,m,l}$ have the same dimension, 
which implies that $S$ is the formal neighborhood of zero 
in ${\mathcal H}_{Y,m,l}-(0,c_0)$. This proves 
part (i) of Theorem \ref{mainth}. 

We now conclude the proof of Theorem \ref{mainth}. 
Let $X'$ be the formal neighborhood
of $M_{c_0}$ in $X$. We have shown that 
the morphism $f: X\to U$ 
lands in ${\mathcal H}_{Y,m,l}$, and that 
$f|_{X'}: X'\to (0,c_0)+S$ is
an isomorphism. This implies that the map $f: X\to {\mathcal H}_{Y,m,l}$ 
is \'etale at $M_{c_0}$. This proves part (ii) 
of Theorem \ref{mainth}, and also implies (iii), since a map 
which is \'etale at one point is dominant.


\begin{thebibliography}{APK}




\bibitem[CBH]{CBH} W. Crawley-Boevey, M. Holland,
{\em Noncommutative deformations of Kleinian singularities},
Duke Math. J. {\bf 92} (1998), no. 3, 605--635. 

\bibitem[CE]{CE} T. Chmutova and P. Etingof, 
On some representations of the rational Cherednik algebra,
math.RT/0303194.

\bibitem[EG]{EG} P. Etingof, V. Ginzburg,
{\em Symplectic reflection algebras, Calogero-Moser space, and 
deformed Harish-Chandra homomorphism}, Invent. Math. {\bf 147} 
(2002), no. 2, 243--348, {\tt math.AG/0011114}.

\bibitem[EO]{EO} P. Etingof, A. Oblomkov,
{\em Quantization, orbifold cohomology, and Cherednik algebras},
preprint, {\tt math.QA/0311005}.

\bibitem[GG]{GG} W. L. Gan, V. Ginzburg 
{\em Deformed preprojective algebras and symplectic reflection algebras for wreath products},preprint, {\tt math.QA/0401038}

\bibitem[GS]{GS} I. Gordon, S.P. Smith
{\em Representations of symplectic reflection algebras and resolutions of deformations of symplectic quotient singularities},
preprint, {\tt math.RT/0310187}.

\bibitem[VB1]{VB1}M. Van Den Bergh
{\em A relation between Hochschild homology and cohomology for Gorenstein rings},Proc. Amer. Math. Soc. {\bf126} (1998), no. 5, 1345--1348. 

\bibitem[VB2]{VB2}M. Van Den Bergh
{\em Erratum to "`A relation between Hochschild homology and cohomology for Gorenstein rings"'},Proc. Amer. Math. Soc. {\bf130} (2000), no. 9, 2809--2810. 



\end{thebibliography}
\end{document}